\documentclass[12pt]{article}
\newcommand{\inst}[1]{$^{#1}$}

\newcommand\email[1]{\small E-mail: #1}

\usepackage{graphicx}
\usepackage{hyperref}
\usepackage{amsmath}
\usepackage{amsfonts}
\usepackage{amssymb}
\usepackage{color}

\usepackage{color}

\begin{document}

\title{Oscillating mushrooms: adiabatic theory\\
for a non-ergodic system%
\footnote{to appear in  J. Phys. A (2014)}}

\author{V. Gelfreich\inst{1},
V. Rom-Kedar\inst{2},
D. Turaev\inst{3}\\
\inst{1}
Mathematics Institute, University of Warwick\\
\email{v.gelfreich@warwick.ac.uk}
\\
\inst{2}
Department of Mathematics, The Weizmann Institute of Science\\
\email{vered.rom-kedar@weizmann.ac.il}
\\
\inst{3}
Department of Mathematics, Imperial College, London\\
\email{dturaev@imperial.ac.uk}
}
\date{}

\maketitle

\begin{abstract}
Can elliptic islands contribute to sustained energy growth as parameters of a Hamiltonian system slowly
vary with time? In this paper we show that a mushroom billiard with a periodically oscillating
boundary accelerates the particle inside it exponentially fast. We provide an
estimate for the rate of acceleration. Our numerical experiments
corroborate the theory. We suggest that a similar mechanism applies to general systems with mixed phase space.
\end{abstract}

\section{Introduction}

Consider a particle which moves freely inside a bounded domain (a billiard) and reflects elastically from
the domain's boundary. We assume that the boundary of the billiard changes with time and restores
its shape periodically. Elastic collisions with the moving boundary cause changes in
 the particle kinetic energy.
If the average energy gain over multiple collisions is positive, the particle accelerates.
This process is often called ``Fermi acceleration'' since it resembles the mechanism of cosmic particles acceleration
by reflecting from magnetic mirrors proposed by Fermi in 1949 \cite{Fermi}.
In  recent years the problem  has attracted a lot of attention (see e.g.
\cite{LR,Dolgopyat2008,LenzPDS2010,STR2010,BR2011,B2014,GRT2012,Shah2011} and references therein).

When the particle accelerates indefinitively, the motion of the boundary eventually becomes 
slow compared to the particle motion. So we assume that 
the separation of  time scales is already present in the billiard.
From the physical point of view, it is important that
the motion of the billiard boundary is not affected by the collisions with the particle, 
i.e. the boundary can be considered as a wall of infinite mass. 
Thus, the Fermi acceleration can be viewed as a process of energy transfer
from a  slow heavy object (in our case, the billiard boundary) 
to fast light particles. In this respect,
the crucial question is whether the  energy transfer is possible at all, and if so, how effective it is, i.e.
what is the rate of the particle energy growth and how it distributes among an ensemble of particles.

We note that there are two alternative approaches to measurements of rates for the Fermi acceleration
in billiards.
Some authors (e.g. \cite{IN03,dCSL2006PRE,LenzPDS2010,LDS08,BR2011})
study the growth of the energy as a function of the collision number $n$. It is easy to see \cite{GRT2012}
that the particle speed can grow at most linearly in $n$. We prefer an alternative point
of view, where the billiard is considered as a flow and we study the energy growth with time $t$.
When a particle accelerates the rate of collisions increases, and the resulting energy
growth rate may be either polynomial or exponential \cite{GT2008,GRST2011,LBPDLS2011}.

It turns out that the acceleration rates  strongly depend on the
shape of the billiard or, more precisely, on the dynamics of 
the corresponding ``frozen'' billiards. At any given moment of time $t$, let us stop the slow motion of
the billiard boundary. 
The dynamics within the frozen domain consist of 
free inertial motion inside the domain and elastic reflections from its static
boundary. This standard billiard dynamical system  is completely
determined by the shape of the frozen domain. Thus, one considers the
family of static billiards parameterised by the frozen time $t$.
When all ``frozen'' billiards in this family are integrable (e.g.
have a rectangular\footnote{The rectangular 
case reduces to the original one-dimensional Fermi-Ulam model \cite{Ulam1961}
in which no acceleration occurs when the boundary oscillations are smooth in time \cite{Pustylnikov1987,Pustylnikov1995}.},
circular or elliptic shape), 
the papers \cite{IN03,dCSL2006PRE,LenzPDS2010,LDS08}
report either no Fermi acceleration or a  slow one.
On the other hand,  if all the frozen billiards are chaotic 
(e.g. Sinai billiards, Bunimovich stadium, etc.),
then the Fermi acceleration is usually present \cite{LRA,LR,KLS2007,B2014a}.

This statement, that chaotic frozen dynamics typically lead to acceleration, is known as the Loskutov-Ryabov-Akinshin
conjecture \cite{LRA} (see also \cite{KLS2007}). In a sense, this conjecture was proven in \cite{GT2008}:
exponentially accelerating trajectories exist in a billiard of periodically oscillating shape provided
every static billiard in the corresponding frozen family has a non-trivial hyperbolic invariant set.
The proof uses the existence of a Smale horseshoe structure in the frozen billiards and, 
as a result, the exponentially fast energy growth is established for initial conditions 
from a set of Lebesgue measure zero only.
Such orbits are difficult to observe numerically and, indeed,  only power-law energy growth
was reported in numerous papers (e.g. \cite{LRA,LR,KLS2007}).
In order to model this moderate growth of the averaged energy
seen in numerical experiments,
 a stochastic differential equation was proposed and  verified numerically in \cite{GRT2012}. According to
this model the ensemble-averaged Fermi acceleration is quadratic in time, provided 
all frozen billiards are ergodic and mixing.
The slow pace of the acceleration is caused by the long-time preservation of the 
Anosov-Kasuga
adiabatic invariant for a large set of initial conditions \cite{Anosov1960,Kasuga1961}. In the billiard context,
the Anosov-Kasuga invariant is equal to $EV^{2/d}$, where $E$ is the particle's energy, $V$ is the volume of the billiard and $d$ is the dimension of the billiard domain\footnote{The 
Anosov-Kasuga theorem is proven only for smooth dynamical systems. 
Billiard dynamical systems have singularities that correspond to corners
and to orbits tangent to the billiard boundary, so the application of the Anosov-Kasuga theory is not formally justified. However,
numerical experiments suggest that this theory is still valid for slowly varying billiards \cite{BOG1987a,GRST2011,GRT2012}}.
For slowly oscillating ergodic billiards, 
the conservation of the Anosov-Kasuga invariant coincides
with the classical thermodynamics adiabatic law for an ideal gas.
In particular, the ergodic adiabatic theory predicts that $\frac{E(T)}{E(0)}\approx 1$ for a large set of initial conditions
if $T$ is the time-period of oscillations of the billiard shape. In this way the ergodic adiabatic theory prohibits
fast acceleration in periodically perturbed ergodic billiards.

A much faster, in fact exponential in time,  
ensemble-averaged Fermi acceleration 
is possible when  
a fraction of the frozen billiard family has several ergodic components \cite{STR2010,GRST2011}. 
Moreover, in this case most  initial conditions experience
exponential energy growth. Namely, transitions between different ergodic 
components of the frozen billiard family lead to a substantial increase in energy transfer to the moving particles.

In this paper we further explore the accelerating effect of the violation of ergodicity:
We consider an example of a slow-fast system whose fast subsystem has a chaotic set coexisting with an elliptic island filled by invariant tori.
We show that  the slow changes in the billiard
shape lead to transitions of the fast variables between the chaotic
and the elliptic zone of the frozen system. We suggest that these transitions  
break the Anosov-Kasuga adiabatic conservation law (which can be related to the Boltzmann entropy
of the system \cite{GRT2012}), cause a systematic increase 
of the entropy over each period of the slow oscillation, 
and thus  lead to a steady (exponential in time) increase in the energy.

We believe that this general principle of induced exponential acceleration 
should be applicable to a wide class of slow-fast Hamiltonian systems whose fast
subspace contains coexisting chaotic and elliptic components \cite{T14}. The recent
preprint \cite{B2014}
provides analytical and numerical arguments  supporting  this
conjecture for billiards with mixed phase space, whereas \cite{PT13}
provides numerical evidence for this phenomenon in a smooth system with mixed phase space.

We stress that this conjecture is of great importance as Hamiltonian 
dynamical systems are rarely ergodic on every energy level. In fact, it is widely believed
that the majority of Hamiltonian systems have mixed phase space.

Here we provide a detailed analysis of this process
%
for a special type of a planar billiard with periodically moving boundaries,
an oscillating Bunimovich mushroom. The frozen billiard shape is shown in Fig.~\ref{Fig:mushroom}. This shape was invented in \cite{Bunimovich2001}. The corresponding billiard has a phase space that is sharply divided into
a single elliptic and a single chaotic component, each one of positive measure \cite{Bunimovich2001,Bunimovich2008}.
One can easily find explicit expressions for the volumes of the regular and chaotic components.
Using this data we propose analytical expressions for the energy distribution
after one cycle of the billiard boundary oscillation and predict the energy growth rate.
We show that provided the billiard boundary moves along a non-trivial loop
in the space of the billiard parameters, a particle inside the billiard accelerates exponentially fast.

Essentially, we show that the process of energy growth over many cycles of boundary oscillations can be modelled
by a geometrical Brownian motion, where the particle energy after each cycle is multiplied by an independent random variable.
We derive an expression for the expectation 
of the logarithm of this random factor 
in terms of the volumes of the chaotic and elliptic
components, and show that this expectation 
is non-negative (and typically strictly positive, which immediately 
implies the exponential energy
growth). Our numerical experiments confirm the predicted growth rate 
with good precision.


The paper has the following structure.
In section 2 we study the properties of the Bunimovich mushroom \cite{Bunimovich2001} 
and find explicit expressions for the volumes of the regular and chaotic components.
Section 3 contains our main theoretical results.
In particular, we derive an adiabatic theory in the presence of particle flux, 
which is used to  calculate the change in the particle's energy
while it stays in the chaotic or regular component.
Then we  find the probability of capture into the regular component. 
Finally, we calculate the averaged growth rate of 
energy over a period of the mushroom oscillations, and show
that this rate is non-negative. 
Moreover, it is strictly positive for generic oscillations. 
Thus the average energy increases exponentially fast.
In section 4 we present several numerical experiments that confirm our prediction
for the energy growth rate, probabilities of capture into the elliptic island 
and distributions of the energy.
 In section 5 we summarise the work and discuss possible extensions for our theory.

\section{Frozen mushroom\label{sec:frozen}}

\begin{figure}
\centerline{
\includegraphics[height=4cm]{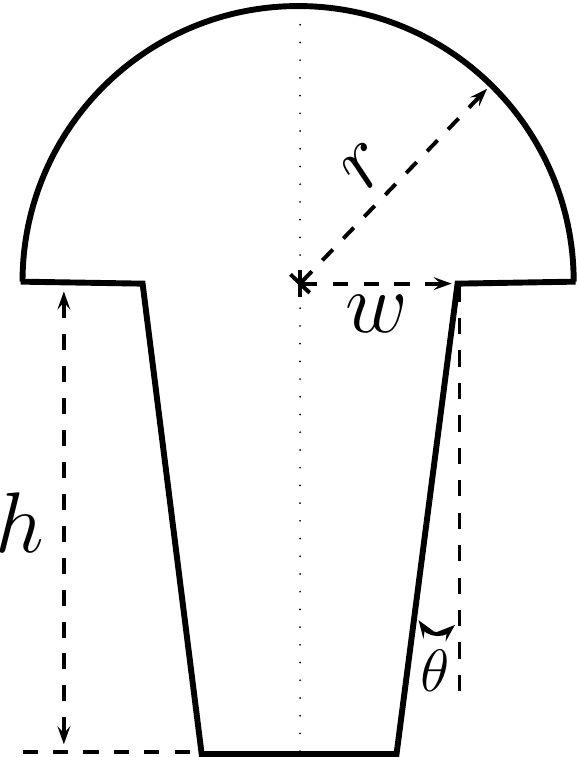}\kern2cm
\includegraphics[height=4cm]{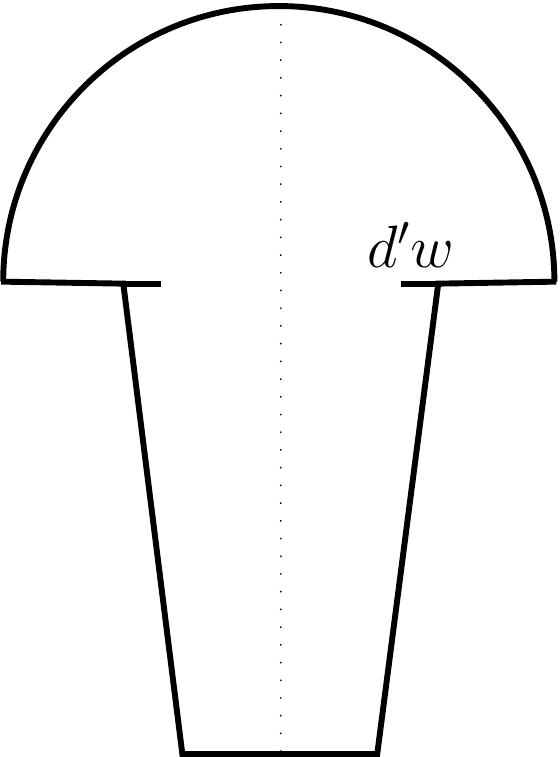}}
\caption{Bunimovich mushrooms}
\label{Fig:mushroom}
\end{figure}

A Bunimovich mushroom consists of a semi-disk and a stem~\cite{Bunimovich2001}
(see Figure~\ref{Fig:mushroom}).
A particle moves freely inside the mushroom $D$
and reflects elastically from its boundary. The particle can go from the cap to the stem and
back through a hole. The dynamics of the particle are defined by the Hamiltonian
$$
H(\mathbf p,\mathbf x)=\frac{p_x^2+p_y^2}{2}\qquad\mbox{for $\mathbf x\in D$}
$$
where $\mathbf p=(p_x,p_y)$ and $\mathbf x=(x,y)$ are  respectively the momentum and position
of the particle. On the boundary of $D$ the particle velocity is reflected
according to the elastic low: the incidence angle is equal to the reflection angle
and the absolute value of the velocity is not changed.

The phase space of the Bunimovich mushroom consists
of two invariant components of positive Lebesgue measure: a regular component (an elliptic island) filled by invariant tori,
and a chaotic component, which is ergodic and mixing \cite{Bunimovich2001,Bunimovich2008}. The regular component
consists of points trapped in the mushroom cap.

Let $r$ and $w$ denote the radius of the cap and the semi-width of the hole respectively
(see Figure~\ref{Fig:mushroom}~(left)).
When the particle moves inside the cap, the incidence angle $\varphi$ with the semi-circular
part of the boundary remains constant.  So if at some collision
$$
|\sin\varphi|\geq\nu=\frac{w}{r}\,,
$$
the particle cannot reach the hole and consequently never leaves the cap.
Along these trajectories the absolute value of the angular momentum, $|x p_y-yp_x|$, is preserved,
and consequently the dynamics are integrable. The chaotic component is the complement
to the integrable one \cite{Bunimovich2008}.

In this paper we will consider a tilted Bunimovich mushroom in the domain $D=D_{\mathrm{stem}}\cup D_{\mathrm{cap}}$
defined by
\begin{eqnarray*}
D_{\mathrm{stem}}&=&\bigl\{\;|x|\leq  w+
y\tan\theta\quad\  \mbox{and }  -h\leq y\leq 0\;\bigr\}\,,
\\[6pt]
D_{\mathrm{cap}}&=&
\bigl\{\;x^2+y^2\leq r \mbox{ and } y\geq 0\;\bigr\}.
\end{eqnarray*}
We assume that $w\leq r$ and $|\theta|<\frac{\pi}{2}$.
For simplicity of presentation we assume that for positive $\theta$ the cone extends up to $y=h$, namely, that $w \ge h \tan\theta$.
At $\theta=0$ we obtain the original Bunimovich mushroom (which has a family of parabolic orbits formed by horizontal oscillations inside the stem,
whereas at nonzero  $\theta$ the family of parabolic orbits is destroyed). It can be shown by methods of 
\cite{Bunimovich1979,Bunimovich2008}
that the tilted mushroom has, for every $\theta$, the same property as the non-tilted one: the complement to the set of points whose orbits are forever
trapped in the cap is ergodic and mixing.

\medskip

Consider a set $D_0\subset D\subset \mathbb{R}^2$ and let  $A(D_0)$ denote its area.
Let $V(D_0)$ be the phase space volume of the set defined by $H(x,y,p_x,p_y)=\frac{1}{2}$ and
$(x,y)\in D_0$. Obviously,
the volume of all points in the phase space such that $H(x,y,p_x,p_y)=E$ and
$(x,y)\in D_0$ is equal to $2\pi\sqrt{2 E} A(D_0)$, so $V(D_0)=2\pi A(D_0)$.

In particular, we find that
\begin{equation}\label{Eq:Vcap}
V_{cap}=\pi^2r^2
\qquad\mbox{and}\qquad
V_{stem}=2\pi (2wh-h^{2}\tan\theta).
\end{equation}
For the future analysis we will need to find volumes of the integrable and chaotic
components. The elliptic island (i.e., the integrable component) occupies only a part of the cap. Let
$V_{ell}$ be  the phase space volume of the elliptic island at the energy level $H=\frac{1}{2}$.
We claim that
\begin{equation}\label{volisl}
V_{ell}=\delta(\nu)V_{cap}(r)
\end{equation}
where $\nu=\frac{w}{r}$ and
\begin{equation}
\delta(\nu)=2\pi^{-1}\left(\arccos\nu-\nu\sqrt{1-\nu^2}\right).
\end{equation}
Indeed, consider  a point $(x,y)=(\rho\cos\psi,\rho\sin\psi)$ inside the cap with a velocity
vector $(p_{x},p_{y})$. Then define the angle $\phi$
by  $(p_x,p_{y})=(v\cos(\psi-\phi), v\sin(\psi-\phi))$  where $v={\sqrt{2E}}$. The absolute value of angular momentum is 
equal to $v\rho|\sin\phi|$.
Notice  that at the cap boundary the angle $\phi$ coincides with the incidence angle $\varphi$.
The integrable component consists of  trajectories that never cross the stem-cap boundary $\{\psi=0,\rho<w\}\cup \{\psi=\pi,\rho<w\}$,
where $w$ is the half width of the hole. Since the angular momentum is conserved, and it is smaller than $wv$ on this boundary,
the integrable component at $v=1$ is defined by the inequalities:
\begin{equation}
{\mathcal D}_{ell}:=\left\{
\frac{w}{|\sin\phi|}\leq \rho\leq r, \quad 0\leq \psi\leq \pi, \quad 0\leq \phi\leq 2\pi\right\}.\label{eq:integrcomp}
\end{equation}
Introducing $s_1=\rho\cos\varphi$, $s_2=\rho\sin\varphi $, we get
$$
{\mathcal D}_{ell}=
\left\{
s_1^2+ s_2^2\leq r^{2}, \quad |s_2|\geq w, \quad 0\leq \psi\leq \pi
\right\},
$$
so the volume is given by the following integral
$$
V_{ell}=\int_{{\mathcal D}_{ell}} d s_1\, d s_2\, d\psi=4\pi \int_w^r \sqrt{r^{2}-s_2^2}\,ds_2\,.$$
Dividing by $V_{cap}$  and using  (\ref{Eq:Vcap}), we obtain (\ref{volisl}).

The chaotic component is the complement of the regular one, therefore its volume is given by
\begin{equation}\label{vchs}
V_{cha}=V_{cap}+V_{stem}-V_{ell}.
\end{equation}

\section{Adiabatic oscillations and capture into an island}

Now suppose that the parameters of the mushroom  change slowly with time.
The speed of the particle is no longer preserved and we study the time evolution of the
particle's energy. The classical adiabatic theory describes dynamics for a slowly changing
integrable system \cite{A1,AKN,LM}, and the ergodic adiabatic theory studies the case when the
system is ergodic at every energy level for every frozen moment of time~\cite{Anosov1960,Kasuga1961}.
Neither of these theories is directly applicable when the phase space of the frozen system
is a mixture of integrable and chaotic components as the particle can transfer from
one type of motion to another one due to the changes in the system.

\subsection{Dynamics in an oscillating circular billiard\label{Se:circular}}

While the particle stays in the cap, its dynamics can be described by a circular billiard
of the same radius. We assume that the radius is a given smooth function of time, $r(t)$, and the centre
of the circle does not move. Let the particle hit the boundary at a point $P$
and $\varphi$ be the impact angle, i.e., the angle between the pre-collision velocity and
the external normal to the boundary at the point of collision. Let $v_{_{\|}}=v\sin\varphi$ denote
the component of the particle velocity parallel to the boundary, and $v_{_{\bot}}=v\cos\varphi$ be the
normal component of the velocity. The elastic reflection $(v_{_{\|}},v_{_{\bot}})\mapsto (\bar v_{_{\|}},\bar v_{_{\bot}})$
from the moving boundary is given by
$$\bar v_{_{\|}}=v_{_{\|}}, \qquad \bar v_{_{\bot}}=2u(t)-v_{_{\bot}}$$
where $u(t)=\dot r(t)$ is the velocity of the boundary motion
(note that  the relation $\bar v_{_{\bot}}-u(t)=-(v_{_{\bot}}-u(t))$ represents the standard elastic law in the coordinate frame that moves
with the boundary). Then for the outgoing angle $\bar\varphi$ we have
$\displaystyle \tan\bar\varphi=\frac{\sin\varphi}{\displaystyle \cos\varphi -2\frac{u}{v}}$
and consequently
\begin{equation}\label{barph}
\bar\varphi=\varphi+2\frac{u(t)}{v} \sin\varphi+O(v^{-2}).
\end{equation}
We assume that the particle moves much faster than the boundary, i.e. $u\ll v$.
The particle speed after the collision is given by
$$\bar v=\sqrt{\bar v_{_{\|}}^2+\bar v_{_{\bot}}^2}=v\left(1-2\frac{u}{v}\cos\varphi+O(v^{-2})\right).$$
The time interval $\Delta t$ to the next collision is $L \bar v^{-1}$ where $L=|PQ|$ is the distance to the next
collision point $Q$. Since  $\Delta t=O(v^{-1})$, the change in the circle radius is also $O(v^{-1})$,
so $L$ is $O(v^{-1})$-close to the value it takes in the static circular billiard of radius $r(t)$,
i.e. to the length of the chord that makes the angle $\varphi$ to the radius. This gives us
$$\Delta t=\frac{2r\cos\varphi}{v}+O(v^{-2}).$$
The new value of the radius at the moment of the next collision is
\begin{equation}\label{rpr}
r'=r(t+\Delta t)=r(t)\left(1+2\frac{u(t)\cos\varphi}{v}\right)+O(v^{-2}).
\end{equation}
Now, by considering the triangle $PQO$ where $O$ is the centre of the circle, we find
$$\frac{r(t)}{\sin\varphi'}=\frac{r'}{\sin\bar\varphi},$$
where $\varphi'$ is the impact angle at the collision point $Q$.
By (\ref{barph}),(\ref{rpr}), we obtain
\begin{equation}
\label{Eq:cbm}
\varphi'=\varphi+O(v^{-2}).
\end{equation}
It follows immediately that the impact angle $\varphi$ stays approximately constant over time required for at least
$O(v)$ collisions, i.e. $\varphi$ is an adiabatic invariant.

Since the circular billiard keeps being rotationally symmetric even when its radius oscillates,
the angular momentum is preserved:
$$
v'r'\sin\varphi'=vr\sin\varphi.
$$
Taking into account equation (\ref{Eq:cbm}) we conclude that
$$
v'r'=vr(1+O(v^{-2}))\,.
$$
It follows immediately that the product $v(t)r(t)$ stays approximately constant over time required for $O(v)$  collisions.
Taking the square, we see that $E(t)V_{\mathrm{cap}}(t)$ is an adiabatic invariant.

In terms of the billiard flow, the adiabatic invariance of the impact angle $\varphi$
can be expressed as the adiabatic invariance of the angle $\hat\varphi$ defined as
\begin{equation}\label{Eq:sinphi}
\sin\hat\varphi=-\frac{x p_y-y p_x}{\sqrt{2E}r(t)}.
\end{equation}
Indeed, at the moments of collision $\hat\varphi$ coincides with $\varphi$,
and it stays approximately constant between impacts.
In the oscillating mushroom the above equations are valid for segments of trajectories
located entirely inside the mushroom cap.

\subsection{Capturing into the cap}

Suppose that a fast particle moves inside the non-autonomous mushroom.
While the particle moves inside the mushroom cap, its dynamics follow the laws of
the circular billiard but the sign of its angular momentum is reflected each time
the particle hits the bottom of the cap.
The hole of the mushroom cap is characterised by the dimensionless parameter
\begin{equation}\label{Eq:nu}
\nu(t)=\frac{w(t)}{r(t)}.
\end{equation}
Suppose that at a moment $t=t_0$ the particle is exactly at the hole. Then it is
located at a point $(x,0)$ with $|x|<w(t_0)$
and has velocity $(p_x,p_y)$. Let $\sin\varphi_0=|\sin\hat\varphi(t_0)|$ be the absolute
value of the adiabatic invariant defined by equation~(\ref{Eq:sinphi}).
Since  $|p_y|\le\sqrt{2E}$
we conclude that $\sin\varphi_0\le \nu(t_0)$.

If at this moment $p_y>0$, the particle enters the cap and can have multiple consecutive collisions
inside the cap. We have already seen that $|\sin\hat\varphi(t)|$ is
adiabatically invariant while the particle stays inside the cap. 
So if the function $\nu$ decreases below $\sin\varphi_0$
while the particle remains inside the cap, the
particle cannot reach the hole. So it must remain in the cap till $\nu$ returns back to $\sin\varphi_0$.

We see that the particle can be captured into the cap when $\nu$ decreases.
If the particle is captured at $t=t_{in}$, then it is released back into
the chaotic zone around the moment $t=t_{out}$ when $\nu(t_{out})=\nu(t_{in})$
for the first time. Thus, we can introduce the release function $t_r(t)$
\begin{equation}\label{trelease}
t_r(t)=\inf\{\, t':t'>t\ \mbox{and}\ \nu(t')\ge \nu(t)\,\},
\end{equation}
which establishes a connection between the capture and release times as
$t_{out}=t_r(t_{in})$.

\subsection{Adiabatic theory in the presence of particle flux}

The classical ergodic adiabatic theory \cite{Anosov1960,Kasuga1961,BOG1987a}
relies on the analysis of evolution of volumes
in the phase space. This theory relies on two observations
which can be presented in a bit oversimplified form in the following way.
First, the ergodicity implies that time averages can be replaced by space averages
and, consequently, the time evolution of the energy is the same for the majority of initial conditions 
starting on a given energy level. Second, if the evolution of the energy depended on initial energy only,
the dynamics would map an energy level into another energy level. Since the Hamiltonian flow is volume preserving, 
the volume under the energy level would stay constant. In a two-dimensional billiard this volume
is proportional to $EV$, which is indeed the ergodic adiabatic invariant. 

If the system is not ergodic on energy levels, the dynamics may produce 
phase space flux between different ergodic components of the frozen system. 
Then, the above volume preservation argument is  invalid. A new paradigm is thus developed

In the non-autonomous mushroom billiard the flux between the regular and chaotic zones
is governed by the parameter $\nu(t)$ defined by (\ref{Eq:nu}).
If $\dot\nu<0$,
 phase volume ``leaks" from the chaotic zone to the regular one,
whereas  $\dot\nu>0$
leads to the opposite effect.

We analyse this situation by considering  a short time interval,  $[t,t+dt]$, at which the mushroom's shape does not change
noticeably, yet, the particle experiences a large number of collisions with
the billiard boundary.
Suppose that during this  time interval the radius of the hole $w$ has changed by $dw$,
the length of the stem $h$ has changed by $dh$ and
the radius of the cap $r$ has changed by $dr$.

In this discussion we need to distinguish two cases depending on the sign of
\begin{equation}
d\nu=\frac{rdw-wdr}{r^2},\label{eq:dnu}
\end{equation}
which determines the direction of the particle flux between the integrable and chaotic
components. In particular, if $rdw<wdr$, the particle in the integrable component
cannot leave the cap
but the particle in the chaotic zone can be captured into the integrable component.
Let us consider the case of capture, $d\nu<0$, in more details.

We conjecture that if the particle is sufficiently fast and the changes in the billiard's shape
are sufficiently small then, from the statistical point of view, the distribution of the energies
at $t+dt$ depends only on the initial and final billiard shapes, i.e., it depends only on the values of
$dw,dh,dr$ and $d\theta$,
and is independent of the particular form of the evolution of the billiard's shape
in the intermediate moments of time.

Now, in order to separate the process of capturing into the cap from the adiabatic
evolution of the energy, we represent the change of the billiard's shape
as a composition of two steps. At the first step,
we allow the mushroom to take the intermediate shape shown on Figure~\ref{Fig:mushroom} (right):
we make the hole in the cap slightly narrower by inserting two straight-line segments into the hole
that symmetrically extend the cup bottom line.

We use the notation $d'z$ and $d''z$ to label changes of a parameter $z$ during the first and
second of these steps respectively ($z\in\{r,w,h,\theta\}$). 
So, over the time interval $dt$ we get $dz=d'z+d''z$.

The shape of the billiard at the end of Stage 1 is uniquely
defined  by the following requirements: On Stage 1 the particle energy
is not changed but all possible transitions between the chaotic and integrable component
take place during this step, i.e., $d'E=0$ and $d'\nu=d\nu<0$.
On Stage 2 the particle cannot move from the integrable
to the chaotic component or vice versa. This condition is achieved by the requirement $d''\nu=0$.
All changes in the energy are on this step: $dE=d''E$.

On Stage 1 all walls remain static (hence $d'h=0$, $d'\theta=0$, $d'r=0$)
except for the bottom of the cap, which
extends to cover a part of the hole as shown on Figure~\ref{Fig:mushroom} (right).
By the end of this stage the radius of the hole is changed by
$$d'w=rd\nu,$$
which is negative since we assume $d\nu<0$. Moreover,  
 the parameter $\nu$ takes its final value:
$$d'\nu=d\nu.$$
The total volume of the phase space remains constant 
on this stage
but a part of the phase volume is transferred from the chaotic to the regular
component:
\begin{equation}\label{Eqdvs}
d'V_{cha}+d'V_{ell}=0,\qquad d'V_{ell}=V_{cap} \;d\delta.
\end{equation}
The last equality is a consequence of equation (\ref{volisl}) and the requirement $d'r=0$.
It is important to note that on Stage~1 the energy remains constant. Indeed, as the straight-line segments
that are inserted into the hole slide along themselves, and the other parts of the billiard boundary
do not move during the Stage 1, the normal velocity of the boundary is non-zero
only at two points, the end points of these two segments. 
Therefore, the particle energy can change only if
it hits the boundary exactly at one of these points at some moment of time, but this is a probability zero event.

On Stage 2 the billiard is slowly deformed from its intermediate shape to the final one in such a way
that the parameter $\nu=\frac{w}{r}$ remains constant, $d''r=dr$ and $d''\nu=0$. Then equation  (\ref{eq:dnu}) implies 
that $d'w+d''w$ attains its correct final value $dw$. At the same time
 $h,\theta$ are changed to ensure  $d''h=dh$, $d''\theta=d\theta$.
It follows that
on Stage 2 the particle is trapped either in the regular or in the chaotic zone
and cannot transfer from one component to the other due to the conservation
of impact angles in the circular part. 

If the particle is in the regular zone, its dynamics are described by the circular billiard
and its energy changes according to the adiabatic
law $d''(EV_{cap})=0$ as derived in Section \ref{Se:circular}.
Since both $E$ and $V_{cap}$ are constant on Stage 1,
and we get $dE=d''E$, $dV_{cap}=d''V_{cap}$ and, consequently,
\begin{equation}\label{Eq:dEe}
\frac{dE}{E}=-\frac{dV_{cap}}{V_{cap}}.
\end{equation}
Since the flux is absent the chaotic zone
remains invariant for the non-autonomous billiard, and we assume that
the standard ergodic adiabatic theory can be applied to the restriction of our dynamical system onto
its invariant subset. So if the particle is inside the chaotic zone, then the absence of flux allows us to claim
that $d''(EV_{cha})=0$.
Using equation (\ref{Eqdvs}) we get
$$
d''V_{cha}=dV_{cha}-d'V_{cha}=dV_{cha}+d'V_{ell}=dV_{cha}+V_{cap}\;d\delta.
$$
Since $d''E=dE$, we  conclude that for particles in the chaotic zone
\begin{equation}\label{Eq:dE}
\frac{dE}{E}=-\frac{dV_{cha}}{V_{cha}}-V_{cap}\frac{d\delta}{V_{cha}}
\end{equation}
where $V_{cha}$ is given by (\ref{vchs}).
Note that in contrast with the classical ergodic adiabatic theory, the right-hand side
of this equality contains an additional term which
takes care of the phase volume flux from the chaotic zone.

One can check that the same equations describe the evolution of the energy
in the case when $d\nu>0$ (flux from the integrable to chaotic component).
The equations can be derived in a similar way
but Stages 1 and~2 are to be swapped.

\subsection{Probability of capture}

In order to describe the acceleration induced by the capture-release
mechanism we consider an ensemble of non-interacting particles
inside the billiard. Then we can discuss the probabilities for a particle to be captured
inside the mushroom cap.

Suppose that $\dot\nu\le0$ on a time interval $[t_0,t_1]$
so the quotient $\nu(t)=w(t)/r(t)$ is decreasing.
We start with $n_{cha}(t_0)$ particles with initial conditions uniformly distributed inside
the chaotic zone.
We assume that the distribution of the particles in the
chaotic zone remains uniform for all times. 

Let $t\in[t_0,t_1)$.
We  consider the transfer of particles from the chaotic
zone during the time interval $[t,t+dt]$
following the two-steps approximation described in the previous section. On Stage 1, a volume of size
$$
d'V_{ell}=V_{cap}\;d\delta
$$
is transferred from the chaotic zone to the regular component,
and on Stage 2 the number of particles in the chaotic zone
remains unchanged. Since the particles are uniformly distributed
we get
$$
\frac{dn_{cha}}{n_{cha}}=-\frac{V_{cap}\;d\delta}{V_{cha}}.
$$
Let $p_{cha}(t)=n_{cha}(t)/n_{cha}(t_0)$ be the probability of being in the chaotic zone
at the time $t$. On the time interval $[t_0,t_1]$, the function
$p_{cha}$ satisfies the relation
\begin{equation}\label{Eq:dps/ps}
\frac{dp_{cha}}{p_{cha}}=\frac{dn_{cha}}{n_{cha}}=
-\frac{V_{cap}}{V_{cha}}\;d\delta.
\end{equation}
Integrating  equation (\ref{Eq:dps/ps}) we  get
\begin{equation}\label{Eq:ps}
p_{cha}(t)=p_{cha}(t_0)\exp\left(
-\int_{t_0}^t\frac{V_{cap}}{V_{cha}}\;d\delta\right)
\end{equation}
for all $t\in[t_0,t_1]$. We stress that these equations are valid only during the period of capture when
no particles are released from the elliptic zone.

\subsection{Energy growth rate over a cycle}

Suppose that the parameters of the mushroom are periodic functions of time and let $T$ be the period.
During this cycle a particle from the chaotic zone may be captured into the island when $\nu(t)$ decreases.
If the particle is captured at a time $t$, then it is released at the time $t_r(t)$
defined by (\ref{trelease}). Note that according to this definition, if $\nu(t)$ is increasing at some $t$, then
$t_r(t)=t$. We also note that $\nu(t)=\nu(t_r(t))$. Since all particles captured between $t$ and $t_r(t)$
are released back into the chaotic zone by the time $t_r(t)$, we also get $n_{cha}(t)=n_{cha}(t_r(t))$ for all $t$.
We assume that at $t=0$ the parameter $\nu(t)$ takes its maximum, so all 
 particles which were captured inside the chaotic
zone during the cycle are releases back by the beginning of the next cycle.

Let us define a compression factor which describes the change in the relative size of the
chaotic zone during the time of capture:
\begin{equation}\label{Eq:g}
g(t)=\frac{V_{cha}(t)}{V_{cap}(t)}\,\left/\,\frac{V_{cha}(t_r(t))}{V_{cap}(t_r(t))}\right..
\end{equation}

Let $S(t)=\log E(t)$ and $\mathbb S(t)=\mathbb E[S(t)]$. Then
equations (\ref{Eq:dEe}) and (\ref{Eq:dE}) imply that
\begin{equation}\label{Eq:dS}
dS=\left\{
\begin{array}{cl}\displaystyle -\frac{dV_{cha}}{V_{cha}}-
V_{cap} \frac{d\delta}{V_{cha}}&\mbox{inside the chaotic zone,}\\[12pt]
\displaystyle-\frac{dV_{cap}}{V_{cap}}&\mbox{in the island.}
\end{array}\right.
\end{equation}
The total number of particles is $n_0=n_{cha}+n_{ell}$. Consequently,
$$
d\mathbb S=
\frac{n_{cha}}{n_0}\left(
-\frac{dV_{cha}}{V_{cha}}
-
V_{cap} \frac{d\delta}{V_{cha}}
\right)-
\left(1-\frac{n_{cha}}{n_0}\right)
\frac{dV_{cap}}{V_{cap}}
\,,
$$
We can  rewrite this equation in the following form
$$
d\mathbb S
=-
\frac{n_{cha}}{n_0}
d\log\left(\frac{V_{cha}}{V_{cap}}\right)
-
\frac{n_{cha}}{n_0}
\frac{V_{cap}}{V_{cha}}
d\delta
- d\log V_{cap}.
$$
Integrating over a complete cycle and taking into account
that $p_{cha}=n_{cha}/n_0$
(all particles are in the chaotic zone at the beginning) we get
$$
\mathbb S(T)-\mathbb S(0)=\int_0^Td\mathbb S(t)=
-\int_0^T p_{cha}
d\log\left(
\frac{V_{cha}}{V_{cap}}
\right)
- \int_0^T p_{cha}
\frac{V_{cap}}{V_{cha}}\,d\delta.
$$
The first term may be integrated by parts:
$$
I_1=
-\int_0^T p_{cha}
d\log\left(
\frac{V_{cha}}{V_{cap}}\right)
=
\int_0^T \log\left(\frac{V_{cha}}{V_{cap}}\right)
dp_{cha}
=
\int_{\mathrm{capture} }
\log g(t)
dp_{cha}, $$
where we grouped together contributions from capture-release pairs
taking into account that $p_{cha}(t)=p_{cha}(t_r(t))$ and $\delta(t)=\delta(t_r(t))$.
For the second integral we get
\begin{eqnarray*}
I_2
&=&-\int_0^T p_{cha}
\frac{V_{cap}}{V_{cha}}\,d\delta
=-\int_{\mathrm{capture} }
p_{cha}(t)
\left(
\frac{V_{cap}(t)}{V_{cha}(t)}
-\frac{V_{cap}(t_r)}{V_{cha}(t_r)}
\right)\,d\delta(t)
\\&=&-
\int_{\mathrm{capture} }
p_{cha}(t)
(1-g(t))
\frac{V_{cap}(t)}{V_{cha}(t)}
\,d\delta(t)
=
\int_{\mathrm{capture} }
(1-g(t))
dp_{cha}
\end{eqnarray*}
where we used equation (\ref{Eq:dps/ps}) which is valid during the capture process.
Defining $p_{ell}=1-p_{cha}$
we get the formula
\begin{equation}\label{Eq:Erate}
m_1:=
\mathbb E\left[\log\frac{E(T)}{E(0)}\right]=\mathbb S(T)-\mathbb S(0)=\int_{\mathrm{capture}}
(g-1-\log g) dp_{ell}\,.
\end{equation}
Since $p_{ell}(t)$ is a non-decreasing  function of time
during the capture
and $\log g<g-1$ for any $g>0$, $g\ne1$,
we conclude that the energy growth rate $m_1$ is non-negative. Moreover, it is strictly positive if $g(t)\not\equiv1$.
As $p_{ell}$ increases during the capture process, we have shown that
$$
\mathbb E\left[\log\frac{E(T)}{E(0)}\right]>0
$$
for any periodic cycle which is nontrivial, namely,
for any cycle having a non-trivial interval of capture with $g(t)\not\equiv1$ on this interval.

Since $\nu(t)=\nu(t_r(t))$ and $\delta$ is a function of $\nu$ only,
we get $\delta(\nu(t))=\delta(\nu(t_r(t)))$, i.e., at the moments of capture
and release, the regular zone takes the same proportion of the cap's phase space volume.
We conclude that $g(t)=1$ if and only if
$V_{cap}(t)/V_{stem}(t)=V_{cap}(t_r)/V_{stem}(t_r)$ (see (\ref{Eq:g}),(\ref{vchs})).

This observation can be restated in the following way.
The equation 
\(\left(\frac{V_{ell}(t)}{V_{cap}(t)},\frac{V_{cap}(t)}{V(t)}\right)\)
with $t\in[0,T]$
defines a closed curve.
If this curve encloses a non-empty interior, the cycle is nontrivial.
On the other hand, any cycle with an empty interior is trivial.
In particular, if one changes only a single 
parameter of the billiard, the above mechanism does not produce
the exponential acceleration and we expect much slower acceleration rates.

\bigskip

Next, we assume that all particles have the same energy $E(0)$ at the beginning of
the billiard cycle and derive an equation for the energy distribution 
at the end of the cycle. For simplicity we assume that the billiard cycle
contains a single interval of capture, i.e. $\dot\nu$ is negative only on
a single interval of time during one complete cycle of the billiard boundary.
Then each particle can be captured at most once per billiard cycle.
Suppose that a particle
is captured at $t=t_{in}$ and let $E_1(t_{in})$ be its energy
at the end of the cycle. In the adiabatic approximation its energy
can be obtained by integrating equation (\ref{Eq:dS}):
\begin{eqnarray}\nonumber
\log\frac{E_1(t_{in})}{E_0}&=&
\int_{[0,t_{in}]\cup[t_{out},T]}
\left(
-\frac{dV_{cha}}{V_{cha}}
-V_{cap} \frac{d\delta}{V_{cha}}
\right)-
\int_{t_{in}}^{t_{out}}
\frac{dV_{cap}}{V_{cap}}
\\
&=&
-\log g(t_{in})+
\int_{t_{out}}^{t_{in}+T}\left(
-
V_{cap} \frac{d\delta}{V_{cha}}
\right)\label{Eq:lnE/E0}
\end{eqnarray}
where the second line uses the definition of the compression factor (equation (\ref{Eq:g})) and the fact that $\int_0^T d(\log V_{cha})=0$.
Since the probability of capture at $t\in[t_{in},t_{in}+dt_{in}]$
is equal to $-\dot p_{cha}(t_{in})dt_{in}$ where $p_{cha}$ is given by
equation (\ref{Eq:ps}), we obtain the probability distribution of the energy
after the cycle in an implicit form. Later in this paper we will
use these implicit equations to reconstruct the distribution of the energy
for specific examples of the mushroom cycles.

On the other hand, if the particle is not captured over the cycle its
energy is defined in the adiabatic approximation (see (\ref{Eq:dE})):
\begin{equation}
\log\frac{E_1^{nc}}{E_0}=
\int_0^{T}
\left(
-\frac{dV_{cha}}{V_{cha}}
-V_{cap} \frac{d\delta}{V_{cha}}
\right)
=-\int_{0}^{T}
V_{cap} \frac{d\delta}{V_{cha}}.\label{Eq:lnEnc/E0}
\end{equation}
This equation implies that, if the particle stays in the chaotic zone over the whole
cycle, its energy at the end of the cycle does not need to be equal to the initial
energy. This conclusion is in a strong contrast with the ergodic case, where the
adiabatic theory predicts that the energy returns close to its initial value
at the end of a cycle. The  changes in the energy
are determined solely by the correction term in (\ref{Eq:dE}), which takes into account the phase flux
due to the non-ergodicity of the frozen billiards. So the phase flux influences
the evolution of the energy even for the particles which never cross to
the regular zone.

Any closed curve in the space of parameters defines two billiard cycles which correspond
to two different directions of motion along the curve. The  values of the integrals (\ref{Eq:lnEnc/E0}) 
for these two cycles have the same absolute value but opposite signs. So the energy of
the non-captured particles may increase or decrease after the completion of the cycle, but in both
cases the  energy averaged over all initial conditions increases.

In the next section we  check numerically  the prediction of
equation (\ref{Eq:lnEnc/E0}) for several examples of billiard cycles.

\section{Examples of billiard cycles}

The theory developed in the previous sections relies upon several assumptions, which
are difficult to prove analytically. In particular, we assume that the ergodic averaging theory
is applicable to the chaotic component in a system which does not satisfy some of the
assumptions of the original averaging theory. The most dramatic violation here is that there are transitions between
the ergodic components - we propose that formula (\ref{Eq:dE}) replaces the usual adiabatic law. Additionally, the billiard
is not a smooth system (see discussion in \cite{GRT2012}).  Finally, we
also assume here that the distribution of the particle in the chaotic zone of the breathing billiard is
close to the stationary uniform one. So we carry out numerical tests to check the correctness
of the theoretical predictions for the energy growth rate and for the distribution of the energy after a cycle
of the billiard boundary.

\subsection{Fixed cap}
\begin{figure}
\begin{center}
\includegraphics[width=0.3\textwidth]{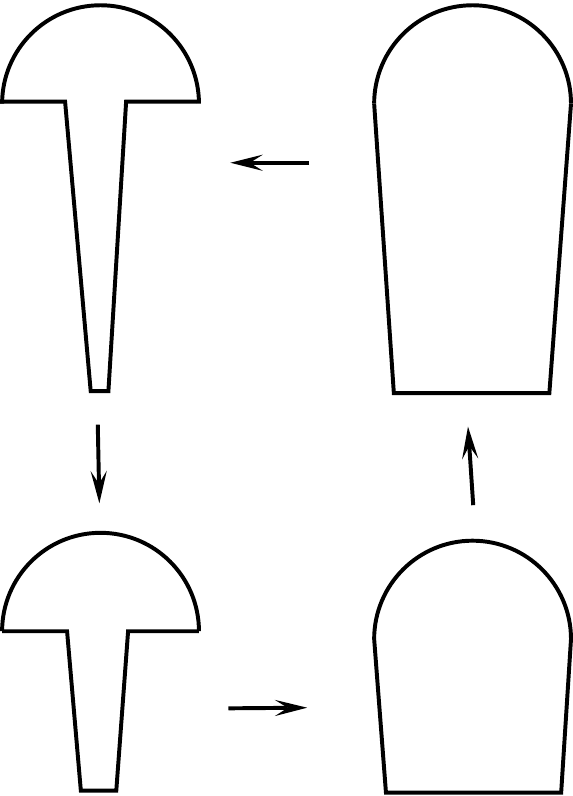}\hskip1.5cm
\includegraphics[width=0.5\textwidth]{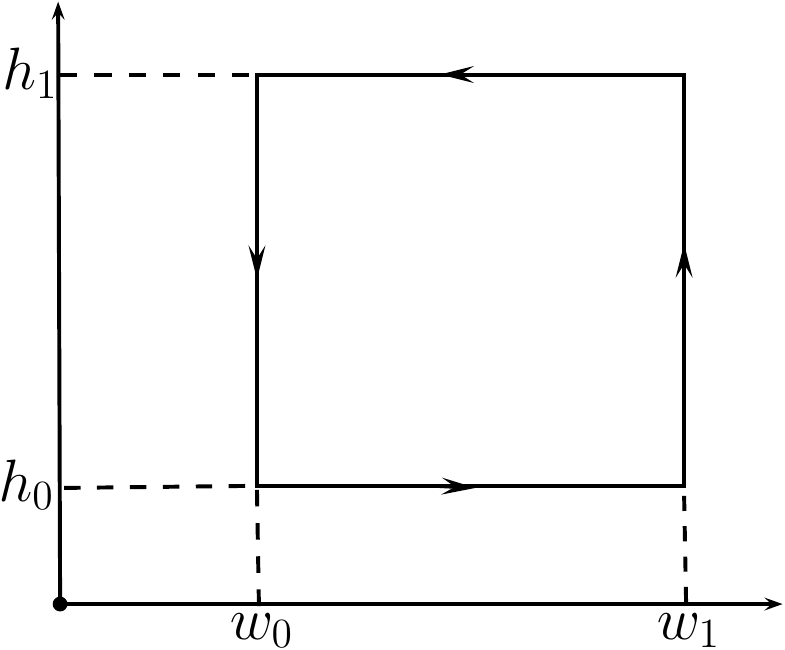}%
\end{center}
\caption{A cycle of the mushroom's oscillations: $h$ and $w$ are the length of the stem and
its radius (half-width) at the cap respectively.}
\label{Fig:cycle}
\end{figure}

Our first example corresponds to the following protocol for the time dependence of the mushroom
parameters: the radius of the gap and the length of the stem follow straight lines
which connect the points
\begin{equation}\label{Eq:cycle}
(w_1,h_{1}
)\rightarrow(w_{0},h_{1})\rightarrow(w_{0},h_{0})\rightarrow(w_1,h_{0})\rightarrow(w_1,h_{1})
\end{equation}
on the plane $(w,h)$ while $\theta$ is fixed and $r=1$ (see Figure~\ref{Fig:cycle}).
For definiteness we assume that $w_0<w_1$ and $h_0<h_1$.
The particle can be trapped inside the mushroom cap
during the first stage of the process, when the hole shrinks,
and then it is released during the third stage.
Since the radius of the cap is fixed, the particle's energy stays constant while
the particle remains in the cap.

The processes of capture and release are determined by the width of the gap:
if the particle is captured into the cap at $t=t_{in}$ then it is released at $t=t_{out}$
such that $w(t_{in})=w(t_{out})$. So
we can rewrite equations (\ref{Eq:Erate}) and (\ref{Eq:g})
\begin{equation}
\mathbb E\left[\log\frac{E_1}{E_0}\right]=
\int_{w_0}^{w_1}
\bigl(g(w)-1-\log g(w)\bigr) dp_{cha}(w)
\end{equation}
where
$$
g(w)=\frac{V_{cha}(w,h_1)}
{V_{cha}(w,h_0)}
\qquad\mbox{and}\qquad
p_{cha}(w)=\exp\left(
-\int_{w}^{w_1}\frac{dV_{ell}(w')}{V_{cha}(w',h_1)}
\right)\,,
$$
and the volumes are defined by (\ref{Eq:Vcap}) and (\ref{volisl}).
In the derivation of the last equality we take into account
that  $V_{cap}$ remains constant and
$V_{ell}=\delta(w)V_{cap}$.

Then $p_{\mathrm{nc}}=p(w_0)$ is the probability of avoiding the capture completely.
If the particle avoids capture,
 equation (\ref{Eq:lnEnc/E0}) implies that the energy after a complete cycle is given by
\begin{equation}\label{Eq:Enoncapt}
\log\frac{E^{\mathrm{nc}}_1}{E_0}=
\int_{w_0}^{w_1}
\left(\frac{1}{V_{cha}(w',h_1)}
-\frac{1}{V_{cha}(w',h_0)}\right)dV_{ell}(w')\,.
\end{equation}
If a particle is captured at the moment when the hole  size is $w$
then its energy at the end of the cycle is described by
(\ref{Eq:lnE/E0}), which takes the form
\begin{equation}
\log\frac{E_1(w)}{E_0}
\label{Eq:Ecapt}
=
\int_{w}^{w_1}
\left(\frac{1}{V_{cha}(w',h_1)}
-\frac{1}{V_{cha}(w',h_0)}\right)dV_{ell}(w')
+\log\frac{V_{cha}(w,h_0)}{V_{cha}(w,h_1)}\,.
\end{equation}
These equations can be used to construct a distribution for $\log\frac{E_1}{E_0}$
at the end of a cycle since they involve integrals of explicitly defined functions.

Another example is obtained when we reverse the protocol (\ref{Eq:cycle})
by following the same path in the space of parameters in the clockwise direction:
\begin{equation}\label{Eq:cyclerev}
(w_1,h_{1}
)\rightarrow(w_{1},h_{0})\rightarrow(w_{0},h_{0})\rightarrow(w_0,h_{1})\rightarrow(w_1,h_{1})
\,.
\end{equation}
The distribution of the energy and the acceleration rate are described by the same equations
but $h_0$ and $h_1$ are swapped.

In the first series of experiments we follow these two protocols by
moving the billiard boundaries with piecewise constant acceleration.
We generate $N=10^5$ initial points uniformly distributed inside the billiard.
All initial conditions have the same energy $E_0=10^9$ and a random direction of the initial velocity.
Then we follow numerically the trajectory for each of the initial conditions during the time required
to complete one cycle of the billiard boundary.

\begin{figure}
\includegraphics[width=0.45\textwidth]{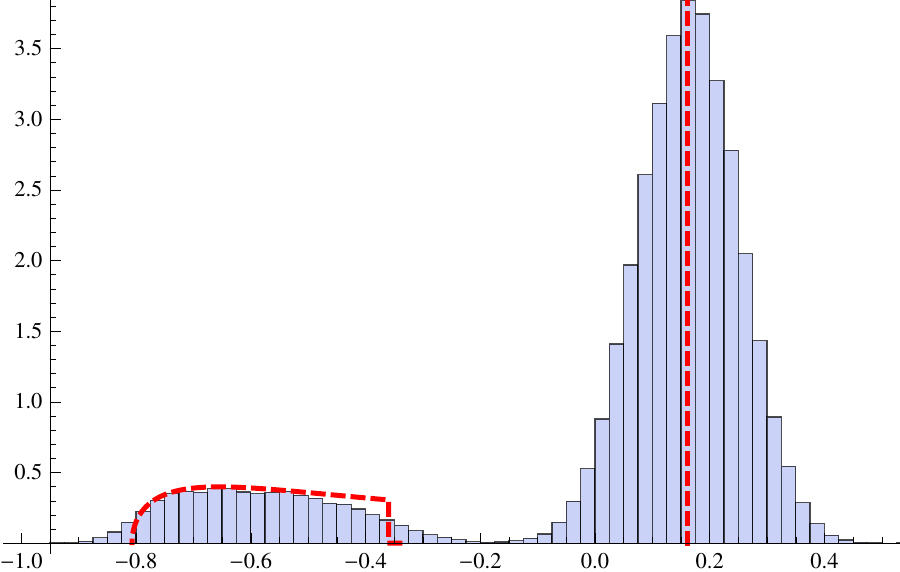}%
\includegraphics[width=0.45\textwidth]{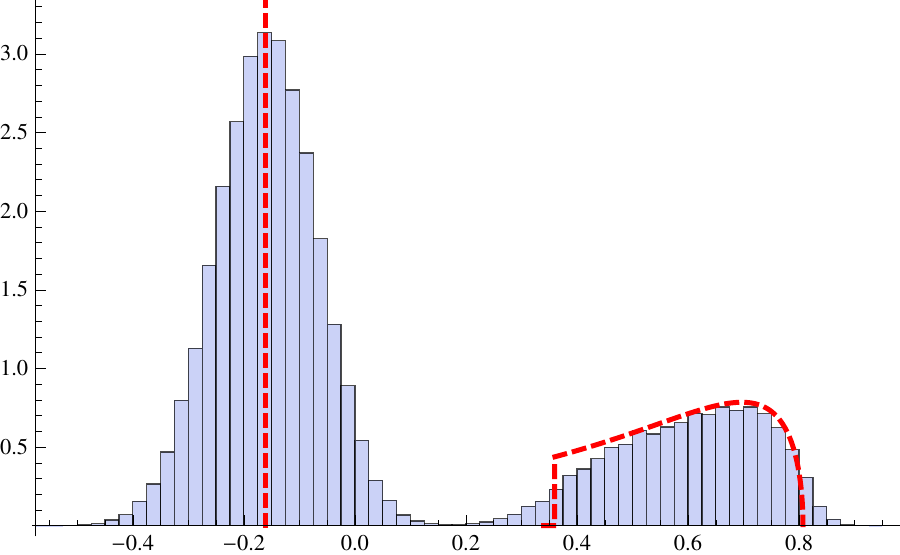}
\caption{Distribution of $\log\frac{E_1}{E_0}$ for anticlockwise (left) and clockwise (right) protocols. The red dashed  line represents
the theoretical prediction. Parameters: $E_0=10^9$, $h_0=2$, $h_1=6$, $w_0=0.3$,
$w_1=1$, $\theta=2.3^\circ$.
\label{Fig:numerics}}
\end{figure}

The distributions of the final energy $E_1$ is
described by histograms which represent relative frequency density for $\log\frac{E_1}{E_0}$.
The histograms are shown on  Figure~\ref{Fig:numerics}
where the dashed red line represents the theoretical predictions for the energy distribution
and the vertical red line marks the position of the theoretically predicted energy
for non-captured particles. In parallel, we mark the capture and release time for each of the
initial condition. In this way we test the accuracy of the theoretical 
prediction for the particle flux into the cap and, in particular, for $p_{\mathrm{nc}}$.

We see that the prediction for the energy growth rate is in good agreement with the numerically
obtained average growth rate:
\begin{itemize}
\item
{\bf Anticlockwise protocol:}
Theoretical prediction: $m_1=\mathbb E[\log\frac{E_1}{E_0}]=0.044926$. 
Non-captured particles: $\log\frac{E_1^{\mathrm{nc}}}{E_0}=0.161205$ with probability $p_{\mathrm{nc}}=0.843472$.
Numerical simulations: average $m_1^*=0.0463\pm0.0027$
for $10^5$ initial conditions.
\item
{\bf Clockwise protocol:} Theoretical prediction: $m_1=\mathbb E[\log\frac{E_1}{E_0}]=0.05252$.
Non-captured particles:
 $\log\frac{E_1^{\mathrm{nc}}}{E_0}=-0.161205$ with probability $p_{\mathrm{nc}}=0.717894$.
 Numerical simulations: average $m_1^*=0.0537\pm0.0034$ for $10^5$ initial conditions.
\end{itemize}
The  histograms of Figure~\ref{Fig:numerics} consist of
two components which correspond to captured and non-captured particles.
The distribution of the energy for the captured particles is apparently in good agreement
with the theory.  The distribution for the non-captured particles looks like Gaussian and its width
scales as $E_0^{1/4}$. This behaviour has been observed for system
where the ergodic averaging theory is applicable \cite{BOG1987a,GRT2012}.
The centre of the distribution is reasonably close to the value predicted by
equation (\ref{Eq:Enoncapt}). We see that the proposed correction to the adiabatic theory (in particular,
equation (\ref{Eq:lnEnc/E0})) is realised. We also
note that numerical experiments show that the relative frequency of escaping the
capture is in excellent agreement with the theoretical prediction given by $p_{\mathrm{nc}}$.
Notice that changing the loop direction leads to changing the role of heating in the stem and cap
--- clockwise motion means that the heating occurs in the cap and the cooling in the stem whereas anticlockwise motion
reverses their role. Yet, as predicted, the overall averaged
growth rate of energy is positive in both cases.

\subsection{Example with moving cap}

In the second set of examples we change the billiard parameters in the following way:
\begin{equation}
\begin{array}{rclrcl}
r(t)&=&r_0+a\sin(t)\,,\qquad&
w(t)&=&r(t)\nu(t)\,,\\[6pt]
h(t)&=&h_0+b\sin( t)\,,&
\nu(t)&=&1-c\sin^2( t)\,.
\end{array}
\end{equation}
In this cycle all parameters of the billiard (except the slope $\theta$)
are changed simultaneously. The capture-release process is determined by
$\nu(t)=\frac{w(t)}{r(t)}$. Since $\nu(t)=\nu(2\pi-t)$ for all $t$ and $\nu$ is monotonically
decreasing on $(0,\pi)$, there is  a 
simple relation between the time of capture and the time of release:
$$
t_r(t)=2\pi-t\qquad \mbox{for $t\in[0,\pi]$}\,.
$$
We find the compression factor from equation (\ref{Eq:g})
and the energy growth rate from (\ref{Eq:Erate}).
Then the distribution of the energy after a complete cycle is
found from  equations (\ref{Eq:lnE/E0}), (\ref{Eq:lnEnc/E0}) and (\ref{Eq:ps}).
For the purpose of plotting the distributions we evaluated these integrals numerically
using the Simpson rule.

\begin{figure}
\begin{center}
\includegraphics[width=0.35\textwidth]{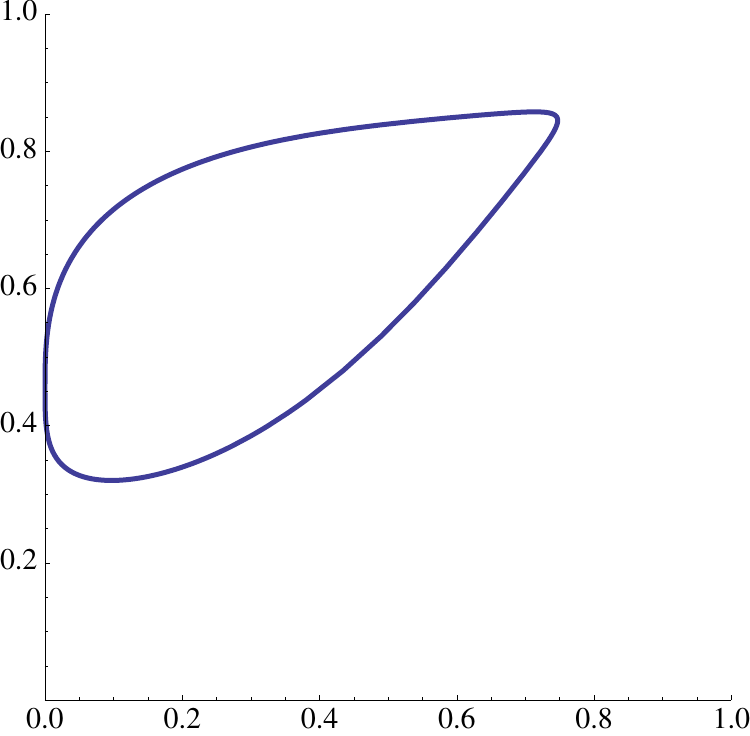}~~\qquad~~~
\includegraphics[width=0.45\textwidth]{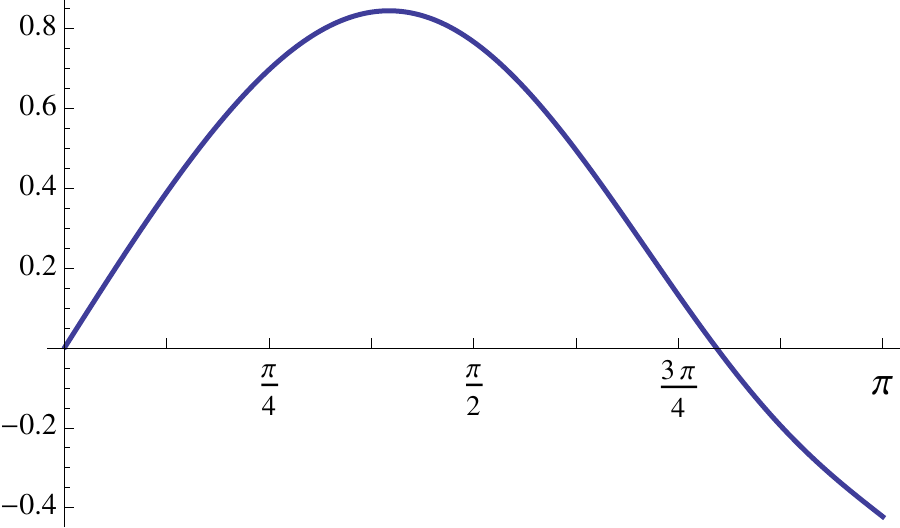}
\end{center}
\caption{The cycle on the plane \(\left(\frac{V_{ell}(t)}{V_{cap}(t)},\frac{V_{cap}(t)}{V(t)}\right)\)
(clockwise direction). Parameters are $r_0=h_0=1$, $a=0.5$, $b=-0.5$, $c=0.8$.
On the right: theoretical prediction for $E_1(t_{in})$ described by equation (\ref{Eq:lnE/E0}).
\label{Fig:protocol}
}
\end{figure}

\begin{figure}
\begin{center}
\includegraphics[width=0.45\textwidth]{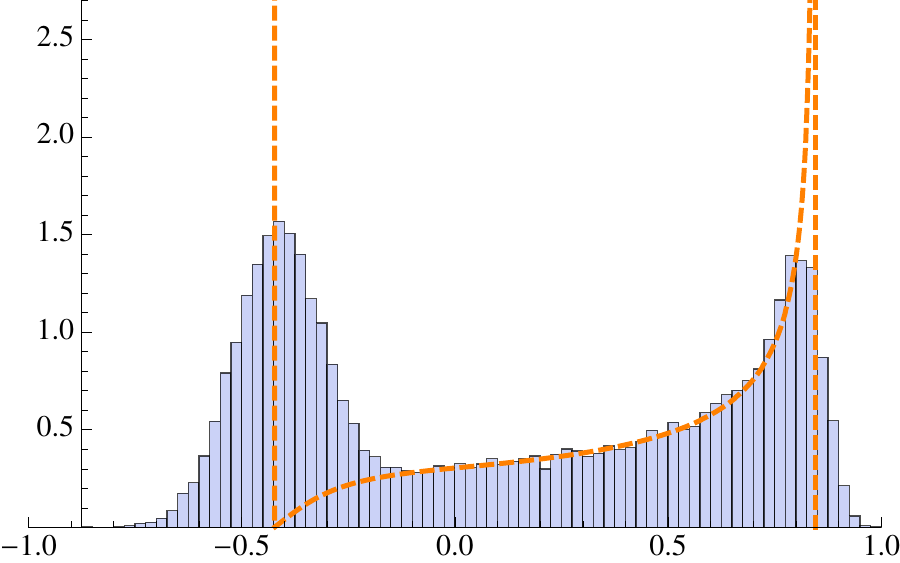}\qquad
\includegraphics[width=0.45\textwidth]{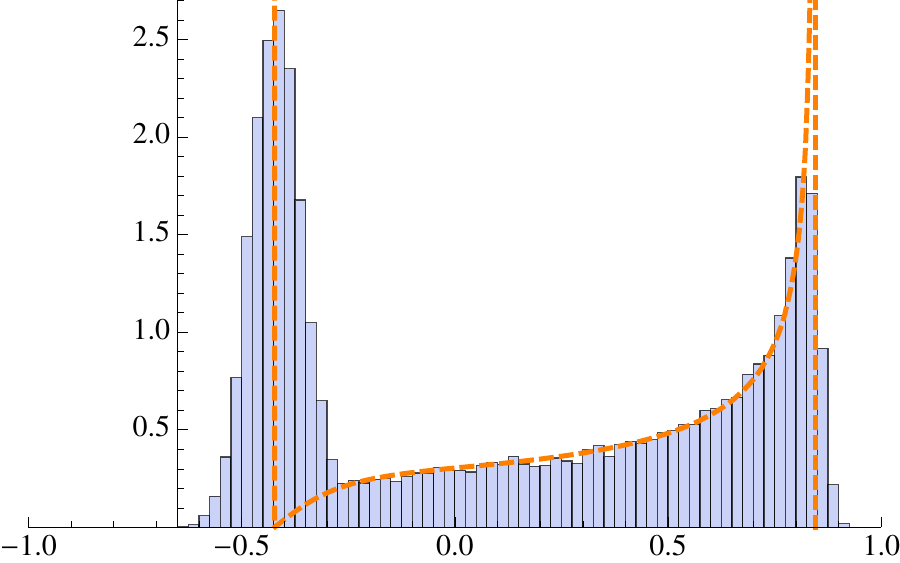}\qquad
\end{center}
\caption{Distributions of $\log \frac{E_1}{E_0}$.
Histograms are constructed on the basis of $25\,000$ initial conditions
uniformly distributed inside the billiard. Initial velocity is taken with random
direction, $E_0=10^6$ (left) and $E_0=10^7$ (right).\label{Fig:E1/E0}}
\end{figure}

In the numerical experiments we use the following values for the parameters:
$$
r_0=h_0=1,\qquad a=-0.5,\qquad b=0.5,\qquad c=0.8\,.
$$
The slope of the billiard stem is kept equal to $\tan\theta=0.1111$
during all experiments of the present section.
The protocol  is illustrated in Figure~\ref{Fig:protocol}.
The right hand side plot
shows the predicted energy after a full cycle for a particle captured at time \(t_{in}\), see equation (\ref{Eq:lnE/E0}).
The positivity of energy gain for most of the captured particles corresponds to the right hump of the distribution shown
in Fig \ref{Fig:E1/E0}. The negative value for \(t_{in}=\pi\)  implies that  
non-captured particles, on
average, loose  energy. A more accurate evaluation of the probabilities suggests
 $\log\frac{E_1^{\mathrm{nc}}}{E_0}=-0.422465$.

We select $N=25\,000$ uniformly distributed initial conditions inside the billiard,
which have the same initial energy and randomly chosen initial directions of velocity.
The distribution of the energy after one cycle is shown on Figure~\ref{Fig:E1/E0}
for two selected values of the initial energy $E_0$. The dashed lines represent
the theoretical distribution. We see that the numerical data are rather close to
the theoretical prediction and, as it should be expected, the agreement is
better for the higher initial energy. The distribution has two modes. The left mode
corresponds to the particles which are not captured into the cap during the cycle.
The position of the left mode is close to the theoretical prediction
given by  $\log\frac{E_1^{\mathrm{nc}}}{E_0}=-0.422465$
while the right mode corresponds to the maximum value of $\log\frac{E_1(t_{\mathrm{in}})}{E_0}$
(compare with Figure~\ref{Fig:protocol}~(right)).

\begin{figure}
\begin{center}
\includegraphics[width=0.45\textwidth]{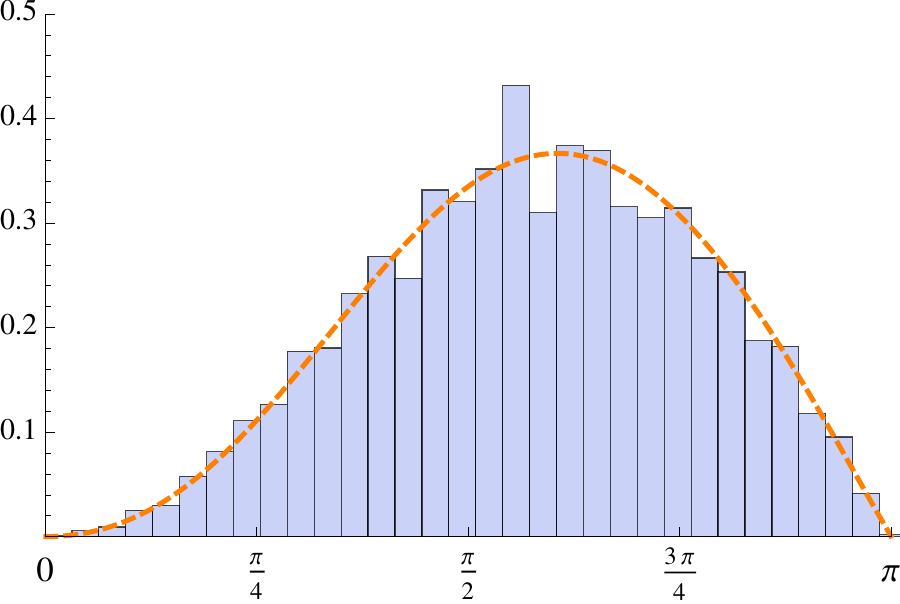}
\end{center}
\caption{Distributions for times of capture\label{Fig:capt}
}
\end{figure}

The theoretical prediction for the growth rate is $m_1=0.122768$. During the numerical
experiment with $E_0=10^7$ we obtain the average growth rate to be
$m_1^*=0.1213\pm0.0033$ $(\pm\sigma_N)$ which is in excellent agreement with the theory.
We also trace the capture and release times for each of the initial conditions.
The distribution of the capture times is illustrated by a histogram shown on Figure~\ref{Fig:capt}.
The theoretical prediction obtained from equation~(\ref{Eq:ps}) is plotted using the dashed line.
In this experiment $38.78\%$ of the particles are not captured in the cap, which is in a good agreement
with the probability of non-capture being $p_{\mathrm{nc}}=38.474\%$.

We conclude that the theory is overall in good agreement with the data from the numerical
experiments, including both the distributions of energy and of capture time.

\begin{figure}
\begin{center}
\includegraphics[width=0.45\textwidth]{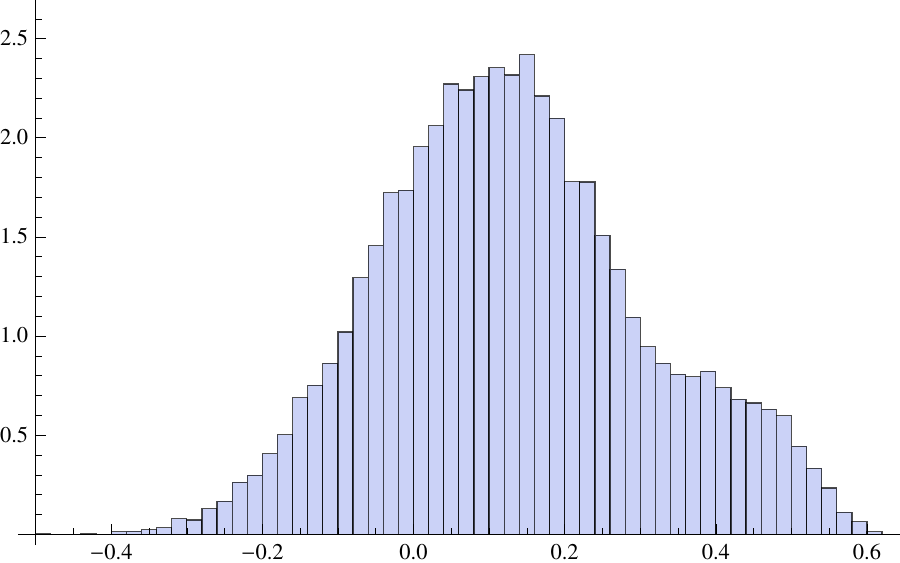}
\qquad
\includegraphics[width=0.45\textwidth]{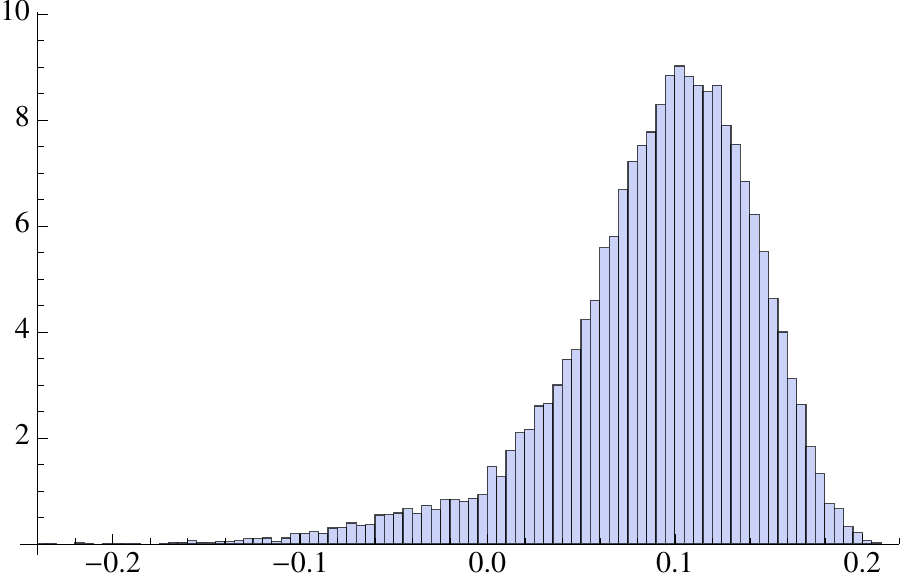}
\end{center}
\caption{Distributions of  $\frac1{10}\log\frac{E_{10}}{E_0}$
and $\frac1{30}\log\frac{E_{30}}{E_0}$.
Histograms are constructed on the basis of $25\,000$ initial conditions
uniformly distributed inside the billiard. Initial velocity is taken with random
direction and $E_0=10^7$.\label{Fig:n=10}}
\end{figure}

We see that after a single cycle the distribution of $\log\frac{E_1}{E_0}$ is quite
far from being Gaussian. The central limit theorem suggest that if the increments
of the logarithm of energy are not correlated over consecutive cycles, then the
distribution of $\frac1n\log\frac{E_n}{E_0}$ should be close to the normal one
centred around $m_1$ for large values of $n$.
Figure~\ref{Fig:n=10} represent the distributions for $n=10$ and $n=30$.
It is clearly seen that while the central part of the distribution is
close to the predicted Gaussian shape, the tails apparently deviate 
from the normal distribution.

\begin{figure}
\begin{center}
\includegraphics[width=0.45\textwidth]{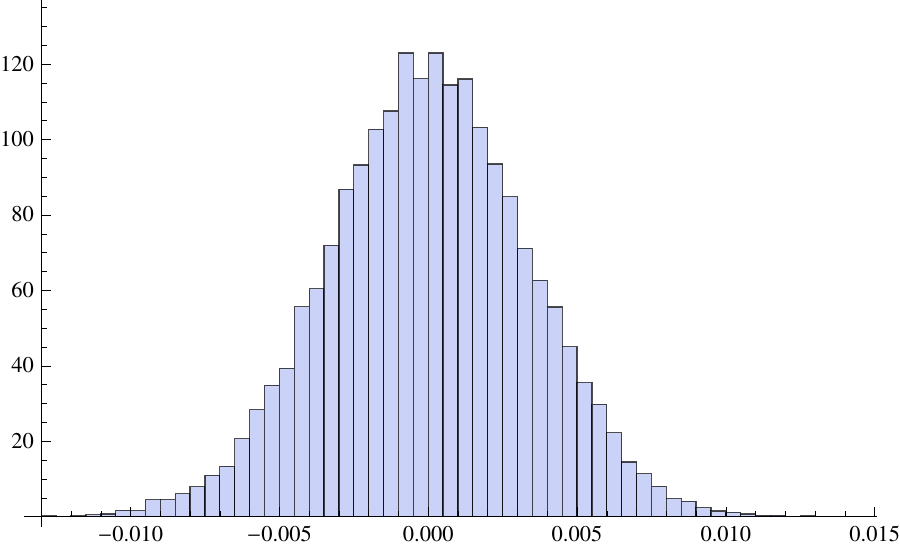}
\end{center}
\caption{Distribution  of $\log\frac{E_1}{E_0}$ for a chaotic billiard.
Histograms are constructed on the basis of $25\,000$ initial conditions
uniformly distributed inside the billiard. Initial velocity is taken with random
direction and $E_0=10^7$.\label{Eq:ergodic}
 $a=0.5,b=-0.5,c=0$, $\tan\theta=0.1111$.
}
\end{figure}

Finally, we stress that the non-ergodicity of the billiard plays the central role
in the creation of the exponential acceleration. In order to illustrate
this difference we consider a billiard cycle with the same parameters as above but
setting $c=0$. While the radius of the cap and the length of the stem
are oscillating as in the previous experiments, 
the width of the stem coincides with the cap diameter and thus the frozen billiard table remains
chaotic at all times. Here
the ergodic adiabatic theory predicts that the energy
should come to its initial value for the majority of the initial conditions.
This is corroborated by a numerical experiment: the distribution of the energy
after one cycle is shown in Figure~\ref{Eq:ergodic}. We see that the final
energy is distributed in a Gaussian-like way with quite small standard deviation:
$|m_1^*|<3\cdot10^{-5}$ and therefore we observe no exponential acceleration on average.

\section{Summary and Discussion}

We propose a mechanism  for achieving an averaged exponential acceleration rate for majority
of initial conditions in a slowly varying system in which the fast dynamics have mixed phase space.
The mechanism is examined by considering an oscillating Bunimovich mushroom, where analytical
predictions are derived and are corroborated numerically. 
In particular, after finding explicit expressions for the  volumes of the regular and chaotic components, 
we derive analytical expressions for the energy distribution
after one cycle and for the exponential energy growth rates. Numerical
experiments support our predictions for both the energy distribution and the growth rate,
and support our claim that the violation of ergodicity is essential for getting exponential acceleration.
We note that  our mechanism does not require precise periodicity
of the process --- it only assumes that the billiard approximately restores its shape and size from time to time.

Our theory involves a generalisation of the ergodic adiabatic theory
which takes into account the flux
between different ergodic components.
We show that the averaged exponential growth rate is described by the following formula
\begin{equation}
m_1:=
\mathbb E\left[\log\frac{E(T)}{E(0)}\right]=\int_{\mathrm{capture}}
(g-1-\log g) dp_{ell}\,.\label{eq:m1conc_}
\end{equation}
where $g$ denotes the compression factor for the phase volumes between the moments of capture into the elliptic island and release
(equation (\ref{Eq:g})) and $p_{ell}$ denotes the probability of the particle being captured in the elliptic component.
This formula is derived under quite general conditions and we hope that it may be
applicable to other systems with mixed phase space.

Averaged exponential acceleration is achieved when  the billiard parameters
change along a non-trivial loop, for which the compression factor $g$ of equation (\ref{eq:m1conc_})
is not identically one.  A non-trivial loop  bounds a non-empty interior when projected on a parameter plane
in which one axis  corresponds to  the phase space volume of the chaotic zone and the other axis corresponds
to the flux between the chaotic and integrable components. An important conclusion is that the exponential acceleration
rate vanishes if the motion of the billiard boundary can be described by periodic
oscillations of a single parameter (as is often done in numerical simulations of Fermi acceleration).
Similarly, it vanishes if one of these two ingredients --- the flux between the ergodic components or
the volume change of the ergodic components --- is missing.
When the exponential rate vanishes, we still expect to observe some slow acceleration, typically quadratic in time (see e.g. \cite{LDS08,LR,LenzPDS2010,LRA}).

We conjecture that the described mechanism, of exponential acceleration due to interior flux and volume changes for different ergodic components of systems
with mixed phase space that are adiabatically deformed, is quite general. However, in contrast with the Bunimovich mushroom, in generic billiards
and in generic smooth systems, 
the  separation of the frozen fast system into  ergodic components that depend continuously on parameters is more problematic.
The existence of chaotic components with positive phase space volume is unknown, and the numerically observed boundary between the seemingly integrable and
chaotic components is often ``sticky'' and fractal.
Nonetheless, we may envision that some rough estimates distinguishing the regular from the chaotic components may be derived (e.g. by calculating Lyapunov exponents),
from which the compression rate $g$ and the capture probability $p_{ell}$ may be found. Then, formula (\ref{eq:m1conc_}) may formally connect these geometrical
features of the frozen system with the energy growth rate in the adiabatically perturbed system.
We should note that the influence of sticky or parabolic orbits on the statistics may be non-trivial. In fact, our initial numerical experiments with
the classical Bunimovich mushroom, which has a rectangular stem and thus a family of parabolic periodic orbits, showed that the energy distribution
in the chaotic zone was different from the one obtained with tilted geometry. The influence of these effects when multiple slow cycles are considered
is yet to be explored.

While we do not anticipate that all the analytical predictions provided here may be carried over to the general case literally, we expect that
the main principle, of achieving exponential acceleration by changing volumes of ergodic components on a non-trivial loop of parameters, is quite general.

\section{Acknowledgements}
This work is partially supported by EPSRC, the Leverhulme Trust, the Royal Society,  the Israel science
foundation (Grant No.  321/12), and the grant RSF 14-41-00044.
\bibliographystyle{unsrt}

\end{document}